\magnification=1200
\hfuzz=3pt
\overfullrule=0mm

\hsize=125mm
\hoffset=4mm


\font\tensymb=msam9
\font\fivesymb=msam5 at 5pt
\font\sevensymb=msam7  at 7pt
\newfam\symbfam
\scriptscriptfont\symbfam=\fivesymb
\textfont\symbfam=\tensymb
\scriptfont\symbfam=\sevensymb

\font\titrefont=cmbx10 at 15pt

\font\refttfont=cmtt10 at 9pt

\font\sc=cmcsc10 \rm


\def\cC{{\cal C}}
\def\Cet{{\cal C_{\hbox{\rm \'et}}}}
\def\cP{{\cal P}}

\def\Hom{{\rm Hom}}
\def\End{{\rm End}}
\def\car{{\rm char}}

\def\ot{\otimes}
\def\eps{\varepsilon}
\def\id{{\rm id}}
\def\card{{\rm card}}
\def\and{\quad\hbox{and}\quad}

\def\bar{\overline}

\def\co{{\rm co-}}

\def\Pr{\noindent {\sc Proof.--- }}

\def\cqfd{ {\sevensymb {\char 3}}}

\hyphenation{commu-ta-tive gene-ra-li-zation theo-re-tical topo-logy
topo-lo-gical
homo-topy homo-topic-ally co-in-va-riant iso-morphism multi-plication
equi-valence equi-valent mor-phism counter-part under-lying tri-viality
tri-vial
affirm-ative invert-ible anti-pode poly-nomials fun-da-mental
non-com-mu-ta-tive
mono-pole con-si-deration pro-po-sition gene-rated}


\null
\noindent
{\sl  
The Legacy of Niels Henrik Abel, The Abel Bicentennial, Oslo, 2002,
O.~A.~Laudal, R.~Piene (eds.), Springer-Verlag 2004, 737--748}

\vskip 40pt

\centerline{\titrefont Quantum principal bundles up to}
\medskip
\centerline{\titrefont homotopy equivalence}

\vskip 30pt

\centerline{\sc Christian Kassel}

\vskip 55pt

\noindent
Hopf-Galois extensions
are known to be the right generalizations of both Galois field extensions and
principal $G$-bundles in the framework of non-commutative associative
algebras.
An abundant literature has been devoted to them by Hopf algebra specialists
(see [Mg], [Sn1], [Sn2] and references therein).
Recently there has been a surge of interest in Hopf-Galois extensions among
mathematicians
and theoretical physicists working in non-commutative geometry \`a la Connes
and \`a la Woronowicz ({\it cf.}\ [BM1], [BM2], [Ha1], [Ha2], [HM], [Mj]).
In their work Hopf-Galois extensions are considered
in the setting of ``quantum group gauge theory".

In this note we deal with Hopf-Galois extensions in the light of topology.
This leads us to simple questions for which we have very few answers, but
which ought to be of interest to those working on Hopf algebras and in
non-commutative
geometry.  We derive these questions from certain
fundamental properties satisfied by topological principal bundles
when we translate them into the setting of Hopf-Galois extensions.
The properties we consider are the following.
\item{(I)} {\it (Functoriality)}
Given a principal $G$-bundle $X\to Y$ and a map $f: Y'\to Y$,
then the pull-back $f^* X  \to Y'$ is a principal $G$-bundle.

\item{(II)}  {\it (Homotopy)} If $f,g : Y' \to Y$ are homotopy equivalent
maps, then $f^*X$
and $g^*X$ are homotopy equivalent bundles.

\item{(III)}  {\it (Triviality)} Any principal $G$-bundle over the point is
trivial.

In order to translate these properties into algebra,
we introduce what we call the homotopy equivalence of Hopf-Galois extensions.
This is the main new concept of this note.
For our definition of homotopy equivalence
we need to restrict to extensions in which the subalgebra of coinvariants
is central.
In other words, the bases of the quantum principal bundles we consider
belong to
classical (commutative) geometry. Nevertheless, we impose no restriction on
the
``structural groups", that is on the Hopf algebras coacting on the quantum
principal
bundles: they may be non-commutative and non-cocommutative,
infinite-dimensional, {\it etc}.

Despite the self-imposed commutativity of the ``base space" the questions
we ask on
Hopf-Galois extensions do not seem easy to answer.
Because they are not the kind of questions usually considered in this subject,
it is our hope that they will lead to new developments on Hopf-Galois
extensions
and encourage research towards their classification over a given base.

We are nevertheless able to answer some of these questions when the Hopf
algebra is one
of  Sweedler or Taft's finite-dimensional Hopf algebras.
This is interesting because the latter are neither commutative,
nor cocommutative, and are sometimes considered as ``toy quantum groups".
In particular, we prove that any cleft extension is homotopy equivalent to
a trivial one.
Cleft extensions form an important class of Hopf-Galois extensions with no
non-trivial counterpart in commutative geometry.
There are some good reasons to view cleft extensions as the natural
non-commutative
generalizations of trivial extensions: like the latter, cleft extensions
can be constructed
from the  ``structural group" and the ``base space"; moreover, in analogy
with the local
triviality of topological bundles, Rumynin~[Ru] showed that
Hopf-Galois extensions over a large class of ``base spaces" become
cleft after (generic) localization. Our above-mentioned result is of a
different, global,
nature; we expect it to hold for any Hopf algebra.

The paper is organized as follows.
In Section~1 we recall the standard terminology related to Hopf-Galois
extensions;
we define what we call quantum principal bundles as well as
their push-forwards along morphisms of commutative algebras (following~[Ru]).
The concept of homotopy equivalence is introduced in Section~2;
we show that up to homotopy equivalence quantum principal bundles
over a polynomial algebra~$C[t]$ are in bijection with quantum principal
bundles over~$C$.
In Section~3 we consider four notions of triviality for a quantum principal
bundle,
and ask whether some of them are equivalent.
In Section~4 we consider the case of Sweedler's and Taft's Hopf algebras
for which
we prove that any cleft quantum principal bundle is homotopy equivalent to
a trivial one.
Section~5 is a short conclusion listing three problems dealing with
Hopf-Galois extensions over the algebra $k[z,z^{-1}]$ of Laurent polynomials,
which we view as the algebraic counterpart of the circle.

We fix a field~$k$.
All constructions will be performed in the category of $k$-vector spaces;
in particular, unadorned tensor symbols refer to the tensor product of
$k$-vector spaces.
All algebras will be assumed to be associative and unital, and all morphisms of
algebras preserve the units.

\bigskip\goodbreak
\noindent
{\bf 1.\ Quantum principal bundles}.
A principal bundle involves a group $G$ (topological group, group scheme)
acting, say on the right, on a space $X$ (topological space, scheme) such
that the canonical
map $X\times G \to X\times_Y X$ induced by $(x,g) \mapsto (x,xg)$ is an
isomorphism
(in the category of spaces under consideration).
Here $Y$ represents some version of the quotient space $X/G$
and $X\times_Y X$ the fibre product.

\medskip\goodbreak
\noindent
1.1.\ {\sc Hopf-Galois extensions}.
In a purely algebraic setting
the group $G$ is replaced by a Hopf algebra $H$ with coproduct
$\Delta : H\to H \ot H$, co\"unit $\eps : H \to k$, and antipode $S : H \to H$.
We assume that the antipode is bijective.
In the sequel we make use of Sweedler's notation: for any $h\in H$ we write
$$\Delta(h) = \sum_{(h)}\, h_{(1)} \ot h_{(2)} .$$

The $G$-space $X$ is replaced by an algebra~$A$ carrying the structure
of an {\it $H$-comodule algebra}. Recall that an algebra~$A$ is an $H$-comodule
algebra if it has a right $H$-comodule structure whose coaction
$\Delta_A : A \to A\ot H$ is a morphism of algebras.
We write
$\Delta_A (a)  = \sum_{(a)}\, a_{(0)} \ot a_{(1)}$ for $a\in A$.

The {\it base} of the $H$-comodule algebra~$A$ is the
subspace~$C = A^{\co H}$ of coinvariants defined by
$$A^{\co H} = \{ a\in A \; | \, \Delta_A(a) = a\ot 1_H \}, $$
where $1_H$ is the unit of~$H$.
It is clear that $A^{\co H}$ is a subalgebra and a subcomodule of~$A$.
We say that $C= A^{\co H}\subset A$ is an {\it $H$-extension}
or that $A$ is an $H$-extension of~$C$.

An $H$-extension $C\subset A$ is said to be {\it Galois}
if the so-called {\it canonical map} $\beta : A\ot_C A \to A\ot H$
defined for $a$, $b\in A$ by
$$a\ot b \mapsto  (a\ot 1_H) \, \Delta_A(b)
= \sum_{(b)}\, ab_{(0)} \ot b_{(1)} \eqno (1)$$
is bijective.
For more details on Hopf-Galois extensions, see [Mg], Chapter~8
or [Sn2], Section~1.

Following Rumynin~[Ru], we consider {\it extensions with central invariants},
name\-ly extensions $C\subset A$ such that $C$ sits in the centre of~$A$
(implying that $C$ is commutative).

The following are the basic objects investigated in this note.

\medskip\goodbreak
\noindent
1.2.\ {\sc Definition}.---
{\it Given a Hopf algebra~$H$ and a commutative algebra~$C$, a
quantum principal $H$-bundle over the base $C$ is a Galois $H$-extension
$C\subset A$ such that $C$ is contained in the centre of~$A$
and $A$ is faithfully flat as a $C$-module. }
\medskip

The faithful flatness condition imposed in the previous definition
is a crucial standard requirement in the subject (see {\it e.\ g.}~[Sn2]).
Observe that our definition of a quantum principal bundle is different
from the one given in [BM2], [Ha1],~[Mj].

\medskip\goodbreak
\noindent
1.3.\ {\sc Functoriality}.
For any morphism of algebras $f: C\to C'$ and any left $C$-module~$A$,
we define the {\it push-forward} $f_*A$ as the left $C'$-module~$C'\ot_C A$.
(This is the algebraic counterpart of a topological pull-back.)
The most natural way to equip $f_*A$ with a multiplication is to define it by
$$(c'_1 \ot a_1) (c'_2 \ot a_2) = c'_1c'_2 \ot a_1a_2$$
for $c'_1$, $c'_2 \in C'$ and $a_1, a_2 \in A$.
Now $c'_1 f(c_1)\ot a_1 = c'_1 \ot c_1a_1$ and $c'_2 f(c_2) \ot a_2 = c'_2
\ot c_2a_2$
in~$C'\ot_C A$ for $c_1, c_2 \in C$. Therefore, in order for the
multiplication
to be well defined, we must have
$$c'_1 f(c_1) c'_2 f(c_2) \ot a_1a_2 = c'_1c'_2 \ot c_1a_1 c_2a_2.$$
This forces us to require $C$ to be central in $A$ and the image $f(C)$ to
be central
in~$C'$. This is why we consider extensions with central invariants.

\medskip\goodbreak
\noindent
1.4.\ {\sc Theorem}.---
{\it
Let $C\subset A$ be a quantum principal $H$-bundle whose base is~$C$.
If $f: C \to C'$ is a morphism of commutative algebras,
then $A' = f_*A$ is a quantum principal $H$-bundle whose base is~$C'$.
}

\medskip
This theorem was proved in [Ru], Theorem~9 under the restriction that
$H$ is finite-dimensional over~$k$.
We give a proof that dispenses with this restriction.
\medskip

\Pr
As we have seen above, $A'$ is an algebra; it contains $C'$ in its centre.
We equip $A'$ with the coaction $\Delta_{A'} = \id_{C'} \ot_C \Delta_A$,
which is well defined because of the condition~$C = A^{\co H}$.
In this way $A'$ becomes an $H$-comodule algebra.

Note that $C' \subset (A')^{\co H}$. This allows us to define a map
$\beta' : A'\ot_{C'} A' \to A' \ot H$ by
$\beta'(a\ot b) = (a\ot 1) \, \Delta_{A'}(b)$ for all $a, b\in A'$.
By assumption, the map $\beta : A\ot _C A\to A\ot H$ defined by~(1)
is bijective. Since
$$A'\ot_{C'} A'  = (C'\ot_C A) \ot_{C'} (C'\ot_C A)\cong  C' \ot_C A \ot_C A$$
and $A'\ot H = C' \ot_C A \ot H$, we have $\beta' = \id_{C'} \ot_C \beta$,
which shows that $\beta'$ is bijective.
By~[Ta], Lemma~4.2, it follows that $C' = (A')^{\co H}$.

It remains to check that $A'$ is faithfully flat as a $C'$-module. This
is a consequence of~[Bo], Proposition~5.
\hfill\cqfd
\medskip

Observe that, if $C\subset A$ is a quantum principal $H$-bundle and $f:C\to
C'$ and
$g:C'\to C''$ are morphisms of commutative algebras, then
$$(g\circ f)_*A = g_*(f_*A).$$
We also have $\id_* A = A$, where $\id:C\to C$ is the identity map.

\bigskip\goodbreak
\noindent
{\bf 2.\ Homotopy equivalence}.
The concept of isomorphism of quantum principal $H$-bundles is clear: two
$H$-extensions
$C\subset A$ and $C'\subset A'$ are isomorphic if there exists an isomorphism
of $H$-comodule algebras $A\cong A'$; such an isomorphism will necessarily
send $C$
onto~$C'$. We'll use the notation $\cong$ for isomorphic quantum principal
bundles.

We now introduce a broader equivalence relation, which is of a global nature.
Let $\Cet$ be the class of {\it finite \'etale morphisms} of commutative
algebras
(for more on \'etale morphisms, see~[Ra]).
We also need the following notation:
if $K$ is a commutative algebra and $K[t]$ is the polynomial algebra
in one variable~$t$ over~$K$, we define $[0], [1] : K[t] \to K$
to be the morphisms of $K$-algebras sending~$t$ to $0$ and~$1$, respectively.

\medskip\goodbreak
\noindent
2.1.\ {\sc Definition}.---
{\it
Given two quantum principal $H$-bundles $A$ and $A'$ over a commutative
base~$C$,
we write $A\sim A'$ if there exist a morphism $i: C\to \bar{C}$ in~$\Cet$ and
a quantum principal $H$-bundle~$B$ with base~$\bar{C}[t]$
such that
$$[0]_* B \cong i_*A \and [1]_* B \cong i_*A'.$$

The homotopy equivalence of quantum principal $H$-bundles over~$C$ is
the equivalence relation generated by~$\sim$.
}
\medskip

Isomorphic quantum principal bundles are clearly  homotopy equivalent.
In Section~4.5 we will show examples of homotopy equivalent
quantum principal bundles that are not isomorphic.
We also remark that by the very definition of homotopy equivalence
two quantum principal bundles $A$ and $A'$ over a base~$C$ such that
$i_*A \cong i_*A'$ for some morphism $i:C\to \bar{C}$ in~$\Cet$ are homotopy
equivalent.

Given a Hopf algebra~$H$ and a commutative algebra $C$, our ultimate goal
is to classify
quantum principal $H$-bundles over~$C$ up to homotopy equivalence.
We denote the set of homotopy equivalence classes of quantum principal
$H$-bundles over~$C$
by~$\cP_H(C)$. The following proposition implies that
the push-forward along a morphism $f: C\to C'$ of commutative algebras induces
a set-theoretic map $f_* : \cP_H(C) \to \cP_H(C')$.

\medskip\goodbreak
\noindent
2.2.\ {\sc Proposition}.---
{\it
Let $f: C\to C'$ be  a morphism of commutative algebras.
If $A$ and $A'$ are homotopy equivalent quantum principal $H$-bundles over~$C$,
then so are $f_* A$ and $f_* A'$.
}

\medskip
\Pr
It is enough to check that, if $A\sim A'$, then $f_* A \sim f_* A'$.
Let $B$ be a quantum principal $H$-bundle with base~$\bar{C}[t]$ ($i: C\to
\bar{C}$
is as in Definition~2.1) such that
$[0]_* B \cong i_*A$ and $[1]_* B \cong i_*A'$.
We extend $f:C\to C'$ to a morphism of algebras
$\bar{f}:\bar{C} \to \bar{C'} = C'\ot_C \bar{C}$. If $i': C' \to \bar{C'}$
is the natural morphism induced from~$i$, we have $\bar{f}\circ i = i'\circ f$.
We next extend $\bar{f}:\bar{C} \to \bar{C'}$ to a morphism
$\bar{f_t} :\bar{C}[t] \to \bar{C}'[t]$ of $k[t]$-algebras.
By Theorem~1.4 $(\bar{f_t})_*B$ is a quantum principal $H$-bundle
over~$\bar{C'}[t]$.
Now, since $\bar{f}(0) = 0$ and $\bar{f}(1)= 1$,
we have $[0]\circ \bar{f_t}  = \bar{f}\circ [0]$ and
$[1]\circ \bar{f_t}  = \bar{f}\circ [1]$. Therefore,
$$[0]_* ((\bar{f_t})_*B) = \bar{f}_*([0]_* B) \cong \bar{f}_*(i_* A) =
i'_*(f_* A)$$
and
$$[1]_* ((\bar{f_t})_*B) = \bar{f}_*([1]_* B) \cong \bar{f}_*(i_* A')
= i'_*(f_* A'). \eqno \hbox{\cqfd}$$
\medskip

We now give a sufficient condition for two morphisms of commutative algebras
$f,g : C\to C'$ to induce the same map $\cP_H(C) \to \cP_H(C')$.
We say that $f$ and $g$ are {\it homotopic}
if there exist a morphism $i': C'\to \bar{C'}$ in~$\Cet$
and a morphism of algebras $\varphi: C\to \bar{C'}[t]$ such that
$[0] \circ \varphi = i'\circ f$ and $[1] \circ \varphi = i'\circ g$
(the morphism $\varphi$ will be called a homotopy between $f$ and~$g$).

\medskip\goodbreak
\noindent
2.3.\ {\sc Proposition}.---
{\it
If $f,g : C\to C'$ are homotopic morphisms of commutative algebras and $A$
is a quantum principal $H$-bundle over~$C$, then the quantum principal
$H$-bundles
$f_*A$ and $g_* A$ are homotopy equivalent: $f_*A \sim g_* A$.
}

\medskip
\Pr
Let $\varphi:  C\to \bar{C'}[t]$ be a homotopy between $f$ and~$g$.
By Theorem~1.4, $A'= \varphi_* A$ is a quantum principal $H$-bundle
over~$\bar{C'}[t]$.
We have
$$[0]_* A'  = [0]_* (\varphi_* A) = ([0] \circ \varphi)_* A = i'_*(f_* A)$$
and
$$[1]_* A' = [1]_* (\varphi_* A) = ([1] \circ \varphi)_* A = i'_*(g_* A) .
\eqno \hbox{\cqfd}$$
\medskip\goodbreak

We give two important applications of Proposition~2.3.
For the first one we let $C = \bigoplus_{i\geq 0} C_i$
be a non-negatively graded commutative algebra
(then $C_0$ is a subalgebra of~$C$).

\medskip\goodbreak
\noindent
2.4.\ {\sc Corollary}.---
{\it The inclusion $\iota: C_0 \to C = \bigoplus_{i\geq 0} C_i$
induces a bijection
$$\cP_H(C_0) \cong \cP_H(C).$$
}

\Pr
Let $\pi: C\to C_0$ be the canonical projection. We have $\pi\circ \iota =
\id_{C_0}$.
We claim the identity map $\id_C$ of~$C$ and $\iota\circ \pi$
are homotopic morphisms. Indeed, let $\varphi: C\to C[t]$ be the morphism
of algebras
defined by $\varphi(x) = t^i x$ if $x \in C_i$;
then $[0] \circ \varphi = i\circ p$ and $[1] \circ \varphi = \id_C$.
It follows from Proposition~2.3 that $A\sim i_*p_* A$
for any quantum principal $H$-bundle $A$ over~$C$.
\hfill\cqfd
\medskip

As a consequence of Corollary~2.4,
the set of homotopy equivalence classes of quantum principal
$H$-bundles over any polynomial ring~$C[x_1, \ldots, x_N]$
with coef\-ficients in a commutative algebra~$C$
is in bijection with the set of homotopy equi\-valence classes of quantum
principal
$H$-bundles over~$C$:
$$\cP_H(C[x_1, \ldots, x_N]) \cong \cP_H(C).$$

\medskip\goodbreak
\noindent
2.5.\ {\sc Corollary}.---
{\it Any morphism $i: C\to \bar{C}$ in~$\Cet$ induces an injection
$\cP_H(C) \to \cP_H(\bar{C})$.
}
\medskip

\Pr
It is enough to check that for any two quantum principal $H$-bundles
$A$ and $A'$ over~$C$ such that $i_* A \sim i_*A'$, we have $A\sim A'$.
Let $\bar{\imath}: \bar{C}\to \bar{C_1}$ be in~$\Cet$ and
$B$ be a quantum principal $H$-bundle with base~$\bar{C_1}[t]$
such that
$$[0]_* B \cong \bar{\imath}_*(i_*A) = (\bar{\imath}\circ i)_*A
\and [1]_* B \cong \bar{\imath}_*(i_*A') = (\bar{\imath}\circ i)_*A'.$$
Since $\bar{\imath}\circ i$ belongs to~$\Cet$ (see~[Ra], Chap.~II,
Proposition~1),
we have~$A\sim A'$.
\hfill\cqfd

\bigskip\goodbreak
\noindent
{\bf 3.\ Triviality}.
In this section we define
four notions of triviality for quantum principal bundles. The first one is the
exact algebraic counterpart of the classical concept of a trivial principal
bundle.
The second one, cleftness,
is proper to the world of non-commutative algebras and has no counterpart
in the classical theory of principal bundles over a group.
The remaining ones, homotopical triviality and homotopical cleftness, are
natural
byproducts of the concept of homotopy equivalence introduced in the
previous section.

\medskip
\noindent
3.1.\ A quantum principal $H$-bundle over a commutative algebra~$C$ is called
{\it trivial} if it is isomorphic to the $H$-comodule algebra $C\ot H$,
where the coaction is given by~$\id_C \ot \Delta$.
(The reader may check that $C\ot H$ is a quantum principal $H$-bundle
over~$C$.)
The push-forward of any trivial quantum principal bundle is trivial.

\medskip\goodbreak
\noindent
3.2.\ An $H$-extension $C\subset A$ is said to be {\it cleft} if there
exists a morphism $\gamma: H \to A$ of right $H$-comodules that is invertible
under convolution in~$\Hom_k(H,A)$, i.e., if there exists another linear map
$\gamma': H \to A$ such that
$$\sum_{(h)}\, \gamma(h_{(1)}) \gamma'(h_{(2)})
= \sum_{(h)}\, \gamma'(h_{(1)}) \gamma(h_{(2)})
= \eps(h) 1_A \eqno (2)$$
for all $h\in H$. The map $\gamma$ is called a {\it cleaving map}.
Note that a cleft extension is called a trivial principal bundle in [BM1],
Example~4.2.

By~[DT1] (see also [Mg], Theorem~7.2.2) an $H$-extension $C\subset A$ is cleft
if and only $A$ is isomorphic to a {\it crossed
product}~$C\,\sharp_{\sigma} H$,
whose underlying vector space is $C\ot H$ and whose multiplication is given
for all $c,d \in C$ and $g,h \in H$ by
$$(c\ot g) (d\ot h)
= \sum_{(g)(h)}\, c(g_{(1)}\cdot d) \,\sigma(g_{(2)}, h_{(1)}) \ot g_{(3)}
h_{(2)}.
\eqno (3)$$
Here $\cdot : H\ot C \to C$ and $\sigma: H\ot H \to C$ are linear maps
defined in terms of the cleaving map $\gamma$ and its inverse $\gamma'$ by
$$h\cdot c = \sum_{(h)}\, \gamma(h_{(1)})\, c \,\gamma'(h_{(2)})  \eqno (4)$$
and
$$\sigma(g,h) = \sum_{(g) (h)}\, \gamma(g_{(1)})\, \gamma(h_{(1)})\,
\gamma'(g_{(2)}h_{(2)})$$
for $g, h\in H$ and $c\in C$.
The coaction of $H$ on $C\, \sharp_{\sigma} H$ is trivial.
A cleft $H$-extension $C\subset A$ is Galois and faithfully flat as a
$C$-module.

For a cleft extension $C \subset C\,\sharp_{\sigma} H$ to be a quantum
principal bundle
in the sense of Definition~1.2, we need $C$ to be central
in~$C\,\sharp_{\sigma} H$. Under
this additional condition, (2) and (4) imply
$h\cdot c = \eps(h) c$ for all $h\in H$ and~$c\in C$.
This means that the ``quasi-action" $\cdot$ is trivial,
which leads to the following simplification of~(3):
$$(c\ot g) (d\ot h)
= \sum_{(g)(h)}\, cd \, \sigma(g_{(1)}, h_{(1)}) \ot g_{(2)} h_{(2)}. $$
Thus, as already observed in~[Ru], Lemma~7,
a cleft quantum principal $H$-bundle is a {\it twisted product}
in the sense of~[Mg], Example~7.1.5.
The push-forward of a cleft quantum principal bundle is cleft by~[Ru],
Theorem~9.

\medskip\goodbreak
\noindent
3.3.\ The third and fourth definitions are as follows:
a  quantum principal $H$-bundle over a commutative algebra $C$ is called
{\it homotopically trivial} (resp.\ {\it homotopically cleft})
if it is homotopy equivalent to a trivial (resp.\ cleft) quantum principal
$H$-bundle over~$C$.
Since push-forwards preserve triviality and cleftness,
Proposition~2.2 has the following consequence.

\medskip\goodbreak
\noindent
3.4.\ {\sc Proposition}.---
{\it
Let $f$ be a morphism of commutative algebras.
If $A$ is a homotopically trivial (resp.\ homotopically cleft)
quantum principal $H$-bundle, then $f_*A$ is homotopically trivial
(resp.\ homotopically cleft).
}
\medskip

The following provides examples of homotopically trivial bundles.

\medskip\goodbreak
\noindent
3.5.\ {\sc Proposition}.---
{\it
If $A$ is a commutative quantum principal $H$-bundle such that the inclusion
$i: C\to A$ belongs to~$\Cet$, then $A$ is homotopically trivial.
}
\medskip

\Pr
Since $A$ is commutative, starting from the quantum principal $H$-bundles
$A$ and
$C\ot H$, the latter being trivial, we may apply Theorem~1.4 and
define the quantum principal $H$-bundles
$i_*A = A\ot_C A$ and $i_*(C\ot H) = A\ot H$ over~$A$. The bundle $i_*(C\ot
H)$ is
trivial. Using the commutativity of~$A$, it is easy to check
that the canonical map $\beta: A\ot_C A \to A\ot H$
is an morphism of $H$-comodule algebras.
(Observe that $H$ is commutative since $\beta$ is an isomorphism.)
Therefore, $i_*A \cong i_*(C\ot H)$ for $i\in \Cet$.
This implies~$A\sim C\ot H$.
\hfill\cqfd
\medskip\goodbreak

We end this section with a few questions.

\medskip\goodbreak
\noindent
3.6.\ {\sc Question}.
{\it Is any cleft quantum principal $H$-bundle over a commutative algebra~$C$
homotopically trivial? }
\medskip

An affirmative answer would give a ``topological" meaning to the algebraic
concept of cleftness and show that homotopical cleftness implies
homotopical triviality.
We may also wonder whether cleftness is closed under homotopy equivalence.

\medskip\goodbreak
\noindent
3.7.\ {\sc Question}.
{\it Let $A$ and $A'$ be homotopy equivalent quantum principal bundles.
Suppose that
$A'$ is cleft. Does it imply that $A$ is cleft?
}
\medskip

An answer in the affirmative to both Questions 3.6 and~3.7 would imply that
the three
notions of cleftness, homotopical cleftness, and homotopical triviality
coincide.
This would give a firm ground to the idea that cleft quantum principal bundles
are the right algebraic counterparts of trivial topological principal bundles
(compare with the question in~[Ru], Section~2.4).

Since a topological principal bundle over a point is trivial,
we ask the following question.

\medskip\goodbreak
\noindent
3.8.\ {\sc Question}.
{\it Is any quantum principal $H$-bundle over the ground field~$k$
homotopically trivial?}
\medskip

According to Schneider~[Sn3], a Galois $H$-extension $A$ of the ground field~$k$ is
cleft in the following cases:

(a) $H$ is finite-dimensional (for a proof, see~[Ru], Corollary 17~(3)),

(b) $H$ is pointed ({\it e.g.}, a group algebra),

(c) $A$ is commutative.

\noindent
Therefore an affirmative answer to Question~3.6 would imply an affirmative one
to Question~3.8 when $H$ satisfies one of the conditions above.

\bigskip\goodbreak
\noindent
{\bf 4.\ The case of Sweedler's and Taft's Hopf algebras}.
The purpose of this section is to answer Questions~3.6 and~3.8 in the
affirmative
for certain small non-commutative non-cocommutative Hopf algebras.

\medskip
\noindent
4.1.\ {\sc Sweedler's four-dimensional Hopf algebra}.
It is the algebra $H_4$ generated by $X$ and $Y$ and the relations
$X^2 = 1$, $Y^2 = 0$, $XY + YX = 0$. The comultiplication of~$H_4$ is given by
$$\Delta(X) = X\ot X \and \Delta(Y) = 1\ot Y + Y\ot X,$$
the co\"unit by $\eps(X) = 1$ and $\eps(Y)=0$,
and the antipode by $S(X) = X$ and $S(Y) = XY$.

For any commutative algebra~$C$ we construct an $H_4$-extension $C\subset A$ as
follows. Let $\alpha$, $\beta$, and $\gamma$ be elements of~$C$ such that
$\alpha$ belongs to the group $C^{\times}$ of invertible elements of~$C$.
We define $A = \Bigl( {\alpha, \beta, \gamma\over C}\Bigr)$
as the $C$-algebra generated by two variables $x$ and $y$ subject to the
relations
$x^2 = \alpha$, $y^2 = \beta$, $xy + yx = \gamma$.
The algebra $A$ is a $H_4$-comodule algebra with $C$-linear coaction
$\Delta_A: A\to A\ot H_4$ defined by
$$\Delta_A(x) = x\ot X \and \Delta_A(y) = 1\ot Y + y\ot X.$$
When $\alpha=1$, $\beta = \gamma = 0$, then
$C\subset \Bigl( {\alpha,\beta,\gamma\over C}\Bigr)$
is trivial in the sense of Section~3.1.

\medskip\goodbreak
\noindent
4.2.\ {\sc Proposition}.---
{\it
The $H_4$-comodule algebra $\Bigl( {\alpha, \beta, \gamma\over C}\Bigr)$ is
a cleft
quantum principal $H_4$-bundle with base~$C$.
Conversely, any cleft quantum principal $H_4$-bundle with base~$C$
is isomorphic to $\Bigl( {\alpha, \beta, \gamma\over C}\Bigr)$ for some
$\alpha\in C^{\times}$ and $\beta$, $\gamma\in C$.
}
\medskip

\Pr
We use Doi and Takeuchi's reformulation of
Masuoka's description of cleft Galois $H_4$-extensions (see [DT2],~[Ms]).
In {\it loc.\ cit.}\ an $H_4$-cleft datum is defined as
a quintuple $(F,D,\alpha, \beta, \gamma)$
where $F$, $D\in \End_k(C)$, $\alpha \in C^{\times}$, $\beta$, $\gamma\in C$
are elements satisfying certain conditions.
Out of any $H_4$-cleft datum
one constructs a cleft $H_4$-extension $C\subset A$,
where $A$ is the free left $C$-module with basis $\{ 1, x, y, xy\}$
equipped with an algebra structure such that
$$x^2 = \alpha, \quad y^2 = \beta, \quad xy + yx = \gamma,$$
$$xc = F(c)x, \quad yc = cy +D(c) x, \quad xyc = F(c)xy + F(D(c)) \alpha$$
for all $c\in C$.
Following~[DT2], we denote $A$ by $\Bigl( {F,D,\alpha,\beta,\gamma\over
C}\Bigr)$
Conversely, any cleft $H_4$-extension $C\subset A$ is isomorphic to
$\Bigl( {F,D,\alpha,\beta,\gamma\over C}\Bigr)$ for some
$H_4$-cleft datum $(F,D,\alpha, \beta, \gamma)$.

Our definition of a quantum principal bundle requires that in addition $C$
be central in~$A$.
This clearly forces $F$ to be the identity map $\id_C$ of~$C$ and $D$ to be
the zero map.
The quintuple $(\id_C,0, \alpha, \beta, \gamma)$ satisfies the requirements
of an $H_4$-cleft datum.
It is clear from the definition of~$\Bigl( {F,D,\alpha,\beta,\gamma\over
C}\Bigr)$
that the algebra associated to the $H_4$-cleft datum $(\id_C,0,
\alpha, \beta, \gamma)$ is $\Bigl( {\alpha, \beta, \gamma\over C}\Bigr)$.
\hfill\cqfd
\medskip

The next result answers Question~3.6 in the affirmative when~$H=H_4$.

\medskip\goodbreak
\noindent
4.3.\ {\sc Theorem}.---
{\it
Any cleft quantum principal $H_4$-bundle over a commutative algebra~$C$
in which $2$ is invertible is homotopically trivial. }
\medskip

\Pr
Let $A$ be a cleft quantum principal $H_4$-bundle over~$C$. By Proposition~4.2
there exist $\alpha\in C^{\times}$ and $\beta$, $\gamma\in C$ such that
$A \cong \Bigl( {\alpha, \beta, \gamma\over C}\Bigr)$.
Let $i : C\to \bar{C}$ be a finite \'etale extension containing a square
root~$s$ of~$\alpha$ (by~[Ra], Chap.~VI, Proposition~1, we may take
$\bar{C} = C[T]/(T^2-\alpha)$).
We have $i_* \Bigl( {\alpha, \beta, \gamma\over C}\Bigr)
= \Bigl( {\alpha, \beta, \gamma\over \bar{C}}\Bigr)$.
By~[DT2], Theorem~2.7 there is an isomorphism
$$\Bigl( {\alpha, \beta, \gamma\over \bar{C}}\Bigr) \cong
\Bigl( {1, \beta, \gamma/s \over \bar{C}}\Bigr).$$
It now suffices to check that any $H_4$-extension of the form
$\Bigl( { 1, \beta, \gamma\over \bar{C}}\Bigr)$ is homotopy equivalent to
the trivial extension $\Bigl( {1,0,0\over \bar{C}}\Bigr)$.
This is proved by considering the cleft quantum principal $H_4$-bundle
$A_t = \Bigl( {1, t\beta, t\gamma \over \bar{C}[t]}\Bigr)$
over the base~$\bar{C}[t]$. We have
$$[1]_* A_t = \Bigl( {1, \beta, \gamma \over \bar{C}}\Bigr)
\and [0]_* A_t = \Bigl( {1,0,0\over \bar{C}}\Bigr). \eqno \hbox{\cqfd}$$
\medskip\goodbreak

This allows us to answer Question~3.8 in the affirmative for~$H_4$.

\medskip
\noindent
4.4.\ {\sc Corollary}.---
{\it  If $\car(k)\neq 2$, then $\card\ \cP_{H_4}(k) = 1$.
}
\medskip

\Pr
By [Ru], Corollary 17~(3), any quantum principal
$H_4$-bundle over the ground field is cleft. We then apply Theorem~4.3.
\hfill\cqfd
\medskip

\medskip\goodbreak
\noindent
4.5.\ {\sc Remark}.
By~[DT2], Corollary~2.8,
$\Bigl( {\alpha, \beta, \gamma\over k}\Bigr)$ is isomorphic to the trivial
quantum principal $H_4$-bundle $\Bigl( {1,0,0\over k}\Bigr)$ if and only if
there exist $s\in k^{\times}$ and $t\in k$ such that
$\alpha = s^2$, $\beta = t^2$, and $\gamma = 2st$.
Therefore, $\Bigl( {1,0,\gamma\over k}\Bigr) \not\cong \Bigl( {1,0,0\over
k}\Bigr)$
for any $\gamma \neq 0$. In view of Corollary~4.4, this shows the existence
of non-isomorphic homotopy equivalent quantum principal bundles.

\medskip\goodbreak
\noindent
4.6.\ {\sc Taft's Hopf algebras.}
Similar results hold for Taft's Hopf algebra $H_{N^2}$, where $N$ is a
fixed integer~$>2$.
Recall that $H_{N^2}$ is the algebra (of dimension~$N^2$)
generated by two generators $X$ and $Y$ subject to the
relations $X^N = 1$, $Y^N = 0$, and $YX = qXY$, where $q$ is a root of
unity of order~$N$.
The comultiplication and the co\"unit are defined by the same formulas as
for~$H_4$.

\medskip\goodbreak
\noindent
4.7.\ {\sc Theorem}.---
{\it
Any cleft quantum principal $H_{N^2}$-bundle over a commutative algebra~$C$
in which $N$ is invertible
is homotopically trivial.
}
\medskip

The proof of Theorem~4.7 relies on [DT2], Theorems 3.2, 3.5,~3.8 and
follows the same lines
as the proof of  Proposition~4.2 and Theorem~4.3.

\medskip
\noindent
4.8.\ {\sc Corollary}.---
{\it  If $N$ is invertible in~$k$, then $\card\ \cP_{H_{N^2}}(k) = 1$.
}

\bigskip\goodbreak
\noindent
{\bf 5.\ Concluding remarks}.
So far we have no example of a non-homotopically trivial quantum principal
bundle.
The search of such an example should be the next step in order to
understand better
the concept of homotopy equivalence introduced in this note.
A good candidate for the non-triviality of~$\cP_H(C)$
may be the algebra $C = k[z,z^{-1}]$ of Laurent polynomials:
it is an algebraic model for the circle, which is the simplest
manifold with non-trivial fundamental group.

\medskip\goodbreak
\noindent
5.1.\ {\sc Problem.}
Find examples of quantum principal bundles over $k[z,z^{-1}]$ that are not
cleft,
not homotopically cleft, or not homotopically trivial.

\medskip\goodbreak
\noindent
5.2.\ {\sc Problem.}
Classify all quantum principal bundles over $k[z,z^{-1}]$ up to isomorphism
(for some interesting Hopf algebras). Find the cleft ones.

\medskip\goodbreak
\noindent
5.3.\ {\sc Problem.}
Determine $\cP_H(k[z,z^{-1}])$ when $H = H_4$ or $H = H_{N^2}$.

\medskip\goodbreak
\noindent
5.4.\ {\sc Remark}.
The quest for affirmative answers to Questions~3.6 and 3.8
may force us to extend the class $\Cet$ to a bigger class $\cC$ of morphisms
(closed under composition as is the case with~$\Cet$).
Proposition~3.5 indicates that $\cC$
should not be too big if we want to keep significant
non-homotopically trivial quantum principal bundles
(see also the example below).

Fix an integer $N\geq 2$ invertible in the ground field~$k$ and
assume that $k$ contains a root of unity $q$ of order~$N$.
Then the extension $k[z,z^{-1}]\subset k[z^{1/N},z^{-1/N}]$
is a quantum principal $H$-bundle for the Hopf algebra $H$
dual to the group algebra of the cyclic group $C_N$
acting on $k[z^{1/N},z^{-1/N}]$ by $\sigma P(z^{1/N}) = P(qz^{1/N})$,
where $\sigma$ is a generator of~$C_N$.
Since this extension is finite and \'etale, it is homotopically trivial
by Proposition~3.5.

\bigskip
\noindent
{\sc Acknowledgement}.
I started working on the problems stated in this note after some
conversations with Vladimir Turaev in early~1999.

\bigskip\bigskip
\centerline{\bf References}
\bigskip

\noindent
[Bo] {\sc N. Bourbaki},
{\it Alg\`ebre commutative}, Herrmann, Paris, 1961.
\smallskip

\noindent
[BM1] {\sc T.~Brzezi\'nski, S. Majid},
{\it Quantum gauge group theory on quantum spaces},
Comm.\ Math.\ Phys.\ 157 (1993), 591--638; Erratum: 167 (1995), 235.
\smallskip

\noindent
[BM2] {\sc T.~Brzezi\'nski, S. Majid},
{\it Line bundles on quantum spheres},
in Particles, fields, and gravitation (\L \'od\'z 1998), AIP Conf.\ Proc., 453,
Amer.\ Inst.\ Phys., Woodbury, New York, 1998, 3--8.
\smallskip

\noindent
[DT1] {\sc Y. Doi, M. Takeuchi},
{\it Cleft comodule algebras for a bialgebra},
Comm.\ Algebra 14 (1986), 808--818.
\smallskip

\noindent
[DT2] {\sc Y. Doi, M. Takeuchi},
{\it Quaternion algebras and Hopf crossed products},
Comm.\ Algebra 23 (1995), 3291--3325.
\smallskip

\noindent
[Ha1] {\sc P. M. Hajac},
{\it Strong connections on quantum principal bundles},
Comm.\ math.\ Phys.\ 182 (1996), 579--617.
\smallskip

\noindent
[Ha2] {\sc P. M. Hajac},
{\it Bundles over quantum sphere and noncommmutative index theorem},
K-Theory 21 (2000), 141--150.
\smallskip

\noindent
[HM] {\sc P. M. Hajac, S. Majid},
{\it Projective module description of the $q$-monopole},
Comm.\ Math.\ Phys.\ 206 (1999), 247--264.
\smallskip

\noindent
[Mj] {\sc S. Majid},
{\it Some remarks on quantum and braided group gauge theory},
in Quantum groups and quantum spaces,
Banach Center Publications, vol.~40,
Institute of Mathematics, Polish Academy of Sciences
Warszawa, 1997, 335--349.
\smallskip

\noindent
[Ms] {\sc A.~Masuoka},
{\it Cleft extensions for a Hopf algebra generated by a nearly primitive
element},
Comm.\ Algebra 22 (1994), 4537--4559.
\smallskip

\noindent
[Mg] {\sc S. Montgomery},
{\it Hopf algebras and their actions on rings},
CBMS Conf.\ Series in Math., vol.~82, Amer.\ Math.\ Soc., Providence, RI, 1993.
\smallskip

\noindent
[Ra] {\sc M.~Raynaud},
{\it Anneaux locaux hens\'eliens},
Lect.\ Notes in Math.\ vol.~169, Springer-Verlag, Berlin, Heidelberg, New
York,~1970.
\smallskip

\noindent
[Ru] {\sc D.~Rumynin},
{\it Hopf-Galois extensions with central invariants and their
geometric properties},
Algebras and Representation Theory 1 (1998), 353--381.
\smallskip

\noindent
[Sn1] {\sc H.-J.\ Schneider},
{\it Principal homogeneous spaces for arbitrary Hopf algebras},
Israel J.~Math.\ 72 (1990), 167--195.
\smallskip

\noindent
[Sn2] {\sc H.-J.\ Schneider},
{\it Hopf Galois extensions, crossed products, and Clifford theory},
in Advances in Hopf algebras (ed.\ J.~Bergen, S.~Montgomery),
Lecture Notes in Pure and Applied Mathematics, vol.~158,
Marcel Dekker, New York, 1994, 267--297.
\smallskip

\noindent
[Sn3] {\sc H.-J.\ Schneider},
Electronic letter to the author (28 April 1999).
\smallskip

\noindent
[Ta] {\sc M. Takeuchi},
{\it On extensions of formal groups by $\mu^A$},
Comm.\ Algebra 5 (1977), 1439--1481.
\smallskip


\vskip 30pt
\line{Institut de Recherche Math\'ematique Avanc\'ee \hfill}
\line{C.N.R.S. - Universit\'e  Louis Pasteur\hfill}
\line{7 rue Ren\'e Descartes \hfill}
\line{67084  Strasbourg Cedex, France \hfill}
\line{E-mail:  {\refttfont kassel@math.u-strasbg.fr}\hfill}
\line{Fax: +33 (0)3 90 24 03 28 \hfill}
\line{http:/\hskip
-2.5pt/www-irma.u-strasbg.fr/\raise-4pt\hbox{\~{}}kassel/\hfill}
\medskip

\bye